\documentstyle{amsart}

\newtheorem{Theorem}{Theorem}[section]

\newtheorem{Lemma}[Theorem]{Lemma}

\newtheorem{Remark}[Theorem]{Remark}

\newtheorem{Definition}[Theorem]{Definition}

\begin{document}

\title{Factorizations of  birational extensions of local rings}

\author{Steven Dale Cutkosky}
\thanks{Research  of Steven Dale Cutkosky  partially supported by NSF}
\author{Hema Srinivasan}

\maketitle

\section{Introduction}

Suppose that $R$ and $S$ are regular local rings such that $S$ dominates $R$ ($R\subset S$ and the maximal ideal $m_S$ of $S$ contracts to the maximal ideal $m_R$ of $R$).

$R\rightarrow S$ is monomial if there exist regular paramaters $x_1,\ldots, x_m$ in $R$, $y_1,\ldots, y_n$ in $S$, an $m\times n$ matrix $A=(a_{ij})$ of rank $m$ whose entries are non negative integers  and units $\delta_i\in S$ such that
$$
x_i=\prod_{j=1}^ny_j^{a_{ij}}\delta_i
$$
for $1\le i\le m$.

Suppose that $P\subset R$ is a regular prime ($R/P$ is a regular local ring) and $0\ne f\in P$.
The regular local ring $R_1=R[\frac{P}{f}]_m$ where $m$ is a maximal ideal of $R[\frac{P}{f}]$ containing $m_R$ is called a monoidal transform of $R$.

Suppose that $V$ is a valuation ring of the quotient field of $S$ which dominates $S$ (and thus dominates $R$). Then given a regular prime $P$ of $R$ (or of $S$) there exists a unique monoidal transform $R_1$ of $R$ (or $S_1$ of $S$) obtained from $P$ such that $V$ dominates $R_1$ (or $V$ dominates $S_1$).

The local monomialization theorem of \cite{C2} and \cite{C4} shows that given an extension $R\rightarrow S\subset V$ as above such that $R,S$ are essentially of finite type over a field $k$ of characteristic zero, there exists a commutative diagram
$$
\begin{array}{llll}
R_1&\rightarrow&S_1&\subset V\\
\uparrow&&\uparrow&\\
R&\rightarrow&S
\end{array}
$$
such that the vertical arrows are products of monoidal transforms and $R_1\rightarrow S_1$ is monomial.

Suppose that we further have that $R\rightarrow S$ is birational (the induced homomorphism of quotient fields is an isomorphism). If $R\rightarrow S$ is monomial and birational, then we can find regular parameters $\overline y_1,\ldots,\overline y_n$ in $S$ such that
$$
x_i=\prod_{j=1}^n\overline y_j^{a_{ij}}
$$
for $1\le i\le n$
(since $B=A^{-1}$ has integral coefficients).

We may now state Abhyankar's local factorization conjecture (page 237 of \cite{Ab}). Suppose that $R\rightarrow S$ is a birational extension of regular local rings of dimension $n\ge 3$ and $V$ is a valuation ring of the quotient field of $S$ such that $V$ dominates $R$. The conjecture is that there exists a commutative diagram 
\begin{equation}\label{eqC1}
\begin{array}{lllll}
&&T&\subset& V\\
&\nearrow&&\nwarrow&\\
R&&\rightarrow&&S
\end{array}
\end{equation}
where the north east and north west arrows are products of monoidal transforms.

It is proven in \cite{Z} and \cite{Ab1} that there is a direct factorization of $R\rightarrow S$ by monoidal transforms if $n=2$. However, 
examples of the failure of a direct factorization of $R\rightarrow S$  by monoidal transforms are given in \cite{Sa} and \cite{Sh} when $n\ge 3$.

The local factorization theorem is proven when $n=3$ (and $R$ is essentially of finite type over a field of characteristic 0) in Theorem A \cite{C1}.

In Theorem 1.9 \cite{C2} it is proven that the local monomialization theorem (Theorem 1.1 \cite{C2}) and ``strong factorization'' of birational toric morphisms of nonsingular toric varieties implies the local factorization theorem in all dimensions (in characteristic zero).

There are two published proofs of ``strong factorization'' of birational toric morphisms, \cite{Mo} and \cite{AMR}. They have both been found to have gaps (as explained in the correction to \cite{AMR}). 

Suppose that $R$ is essentially of finite type over a field. In \cite{C2}, a strong version of local monomialization is used to reduce the proof of local factorization to the following problem,
which is essentially in linear algebra.

We assume that
$R\rightarrow S$ is monomial, with respect to regular parameters $x_1,\ldots, x_n$ in $R$ and $y_1,\ldots, y_n$ in $S$, the value group of $V$ is contained in ${\bf R}$, and if $\nu$ is a valuation of the quotient field of $S$ whose valuation ring is $V$, then
$$
\tau_1=\nu(y_1),\ldots, \tau_n=\nu(y_n)
$$
are rationally independent real numbers.

In this special case, we can assume that $R=k[x_1,\ldots,x_n]_{(x_1,\ldots,x_n)}$ and
$S=k[y_1,\ldots,y_n]_{(y_1,\ldots,y_n)}$ where $k$ is a field. We have expressions 
$x_i=\prod_{j=1}^ny_j^{a_{ij}}$ for $1\le i\le n$. If
$$
f=\sum \alpha_{i_1,\ldots,i_n}y_1^{i_1}\ldots y_n^{i_n}\in k[y_1,\ldots, y_n],
$$
we have $\nu(f)=\text{min}\{i_1\tau_1+\cdots+i_n\tau_n\mid \alpha_{i_1,\ldots,i_n}\ne 0\}$.
We will call the local factorization conjecture in this special case the ``monomial problem''.

When $n=3$, the monomial problem is solved by Christensen \cite{Ch}.  In Theorem 1.6 \cite{C2}, the first author extends this to prove a weaker form of the monomial problem for all $n$. By combining this with the local monomialization theorem of \cite{C2}, it was proved in \cite{C2} that a birational extension $R\rightarrow S$ can be factored by $n-2$ triangles of monodial transforms.

Recently, there has been a proof by Karu \cite{K} of this monomial problem. His proof is geometric in nature.

In this paper, we give a self contained proof of the monomial problem. We solve the problem in the spirit of Christensen's original theorem in dimension 3.  In particular, the problem can be stated completely in the language of linear algebra, and we prove it using only linear algebra. As a result, we give an explicit algorithm for the solution of the monomial problem.
This theorem (Theorem \ref{Theorem1})
 is proven in Section 2 of this paper. The solution to the monomial problem is given in Theorem \ref{Theorem15}.
 
 We show in Theorem \ref{Theorem24} of Section 3 of this paper how the local monomialization theorem, Theorem 1.1 \cite{C2} and Theorem \ref{Theorem1} of this paper
  prove the local  factorization conjecture.  
This provides a complete proof to Theorem 1.9 of \cite{C2}.  

A monoidal transform affects the coefficient matrix $A$ as a column addition. 
The valuation can be understood as a column vector $\vec v$ of positive rational numbers. To preserve the property that the valuation ring dominates the monodial transform of the local ring,
 we allow only those column operations on $A$ that keep both $A$ and $A^{-1}\vec v$ positive. We contruct an algorithm here for
finding a sequence of permissible column additions and  interchanges to be followed by a sequence of permissible  subtractions that results in the identity matrix.

\section{Matrix Factorization}
Suppose that $A=(a_{ij})$ is an $n\times n$ matrix with coefficents which are non-negative integers and $\text{Det}(A)=\pm1$.
Further suppose that
$\vec v=(v_1,\ldots,v_n)^t$ is a $n\times 1$ column vector with coefficients which are positive rationally independent real numbers, and 
$\vec w=(w_1,\ldots,w_n)^t=A^{-1}\vec v$ is a vector with positive coefficients (which are necessarily rationally independent).
$(A,\vec v,\vec w)$ satisfying these conditions will be called an $n$-dimensional triple.

The column addition $C_{ij}$ of $A$ which adds the $j$-th column of $A$ to the $i$-th column is called {\bf permissible} 
for $(A,\vec v,\vec w)$ if $w_j-w_i>0$.
The triple $(A,\vec v,\vec w)$ is transformed under the permissible column addition $C_{ij}$ to the triple  $(A(1),\vec v(1),\vec w(1))$, where 
$A(1)=(a(1)_{ij})$  is obtained from $A$ by adding the $j$-th column of $A$ to the $i$-th column, $\vec v(1)=\vec v$ and
$\vec w(1)=(w(1)_1,\ldots, w(1)_n)^t=A(1)^{-1}\vec v(1)$.
$\vec w(1)$  is obtained from $\vec w$ by subtracting the $i$-th coefficient $w_i$ from the $j$-th coefficient $w_j$ of $\vec w$.

The row subtraction $R_{ji}$ of $A$ which subtracts the $i$-th row of $A$ from the $j$-th row is called {\bf permissible} for $(A,\vec v,\vec w)$ if $a_{jk}\ge a_{ik}$
for $1\le k\le n$. The triple $(A,\vec v,\vec w)$ is transformed under the permissible row subtraction $R_{ji}$ to the triple $(A(1),\vec v(1),\vec w(1))$, where 
$A(1)=(a(1)_{ij})$  is obtained from $A$ by subtracting the $i$-th row of $A$ from the $j$-th row, $\vec v(1)$ is obtained from  $\vec v$
by subtracting the $i$-th coefficient $v_i$ from the $j$-th coefficient $v_j$ and 
$\vec w(1)=(w(1)_1,\ldots, w(1)_n)^t=A(1)^{-1}\vec v(1)$. 
We have that
$\vec w(1)=\vec w$.

The row interchange  $T_{ij}$ of $A$ interchanges the $i$-th and $j$-th rows of $A$. $T_{ij}$ transforms the triple $(A,\vec v,\vec w)$
into the triple $(A(1),\vec v(1),\vec w(1))$, where $A(1)$ is obtained from $A$ by interchanging the $i$-th and $j$-th row,
$\vec w(1)=\vec w$ and $\vec v(1)$ is obtained from $\vec v$ by interchanging the $i$-th and $j$-th row of $\vec v$.

In this section, we prove the following theorem:

\begin{Theorem}\label{Theorem1}  
Suppose that $A=(a_{ij})$ is an $n\times n$ matrix with coefficents which are non-negative integers and $\text{Det}(A)=\pm1$.
 Further suppose that
$\vec v=(v_1,\ldots,v_n)^t$ is a $n\times 1$ vector with coefficients which are positive rationally independent real numbers. 
Then there exists a sequence of permissible column additions and row interchanges
$$
(A,\vec v,\vec w)\rightarrow (A(1),\vec v, \vec w(1)) \rightarrow \cdots \rightarrow (A(s), \vec v,\vec w (s))
$$
followed by a sequence of permissible row subtractions 
$$
(A(s),\vec v, \vec w(s))\rightarrow (A(s+1),\vec v(s+1), \vec w(s)) \rightarrow \cdots \rightarrow (A(t), \vec v(t), \vec w(t))
$$
such that $A(t)$ is the $n\times n$ identity matrix.
\end{Theorem}

We will denote the inverse of a matrix $A$ by $B=(b_{ij})=A^{-1}$. 
If a permissible column addition $C_{ij}$ is  performed by adding the $j$-th column of $A$ to the $i$-th column, 
with a resulting transformation of triples $(A,\vec v, \vec w)\rightarrow (A(1),\vec v, \vec w(1))$, then
$B(1)=(b(1)_{ij})=A(1)^{-1}$ is obtained from $B=A^{-1}$ by subtracting the $i$-th row of $B$ from the $j$-th row,
since $C_{ij}^{-1}=R_{ji}$ and 
$$
B(1)=A(1)^{-1}=(AC_{ij})^{-1}=C_{ij}^{-1}A^{-1}=R_{ji}A^{-1}.
$$
Similarly, if a permissible row subtraction $R_{ji}$   is performed by subtracting the $i$-th row of $A$ from the $j$-th row, 
with a resulting transformation of triples $(A,\vec v, \vec w)\rightarrow (A(1),\vec v(1), \vec w)$,
then
$B(1)=(b(1)_{ij})=A(1)^{-1}$ is obtained from $B=A^{-1}$ by adding the $j$-th column of $B$ to the $i$-th column.

If a permissible row interchange $T_{ij}$ is performed, then $B(1)=A(1)^{-1}$ is obtained from $B$ by interchanging the $i$-th and $j$-th
column.

Given a triple $(A,\vec v,\vec w)$, we define $\beta=\max_k\{|b_{k1}|\}$.  We will write $A=(C_1,\ldots, C_n)$.

To simplify notation, we will denote the inverse of a matrix $A(t)$ by $B(t)=(b_{ij}(t))$, and define $\beta(t)=\max_k\{|b_{k1}(t)|\}$.
We will denote $A(t)=(C_1(t),\ldots,C_n(t))$.

\begin{Remark}\label{Remark2}
Fix $i$ and $j$. Either $C_{ji}$ is permissible or $C_{ij}$ is permissible. If $C_{ij}$ is permissible, then after performing $C_{ij}$ a finite
number of times, $C_{ji}$ becomes permissible.   This is because  $C_{ij}$ decreases
$w_j$ by a positive integral multiple of $w_i$.
\end{Remark}

\begin{Definition}\label{Definition3} A permissible $C_{ij}$ is allowable for the triple $(A,\vec v,\vec w)$
if $b_{i1}$ and $b_{j1}$ are both non-zero and have the same sign.
\end{Definition}

\begin{Definition}\label{Definition7} A permissible $C_{ij}$ is *-allowable for the triple $(A,\vec v,\vec w)$
if either $b_{i1}b_{j1}=0$, or $C_{ij}$ is allowable.
\end{Definition}

\begin{Remark}\label{Remark3}
\begin{enumerate}
\item If we perform a *-allowable $C_{ij}$ on the triple $(A,\vec v,\vec w)$ to get $(A(1),\vec v,\vec w(1))$, then $b_{j1}(1)=b_{j1}-b_{i1}$,
$b_{k1}(1)=b_{k1}$ if $k\ne j$ and thus
$$
\beta(1)=\max_k\{|b_{k1}(1)|\}\le \max_k\{|b_{k1}|\}=\beta.
$$
\item Suppose that we fix $i$ and $j$. Then after a finite sequence of consecutive allowable $C_{ij}$ and $C_{ji}$, both $C_{ij}$ and $C_{ji}$ are not
allowable. If at least one of $b_{i1}$, $b_{j1}$ is nonzero, then after a finite sequence of consecutive *-allowable $C_{ij}$ and $C_{ji}$, both
$C_{ij}$ and $C_{ji}$ are not *-allowable.
\end{enumerate}
\end{Remark}

\begin{pf}(of 2)
If  $b_{i1}$ and $b_{j1}$ are nonzero of the same sign, and  we perform $C_{ij}$ (or $C_{ji}$)
to obtain the new triple $(A(1),\vec v,\vec w(1))$, and $b_{i1}(1)$, $b_{j1}(1)$ have the same sign,
we then obtain that $(|b_{i1}|,|b_{j1}|)> (|b_{i1}(1)|,|b_{j1}(1)|)$ in the Lex order on ${\bf Z}^2$. 

Suppose that $b_{i1}\ne 0$ and $b_{j1}=0$. If  $C_{ij}$ is *-allowable, then after performing $C_{ij}$, we obtain that both $C_{ij}$ and $C_{ji}$
are not *-allowable. If $C_{ji}$ is *-allowable, and we perform $C_{ji}$, then $b_{i1}(1)=b_{i1}$, $b_{j1}(1)=0$. By Remark \ref{Remark2}, we can
only perform $C_{ji}$ a finite number of consecutive times.
\end{pf}

\begin{Lemma}\label{Lemma5} There exists a sequence of allowable column additions
$$
(A,\vec v,\vec w)\rightarrow \rightarrow (A(1), \vec v,\vec w(1))\rightarrow \cdots \rightarrow (A(t),\vec v,\vec w(t))
$$
 such that at most two entries of the first column of $B(t)$ are nonzero.
\end{Lemma}

The proof of this lemma is immediate from Theorem 6.3 \cite{C2}.

\begin{Lemma}\label{Lemma6} There exists a sequence of *-allowable column additions
$$
(A,\vec v,\vec w)\rightarrow  (A(1), \vec v,\vec w(1))\rightarrow \cdots \rightarrow (A(s),\vec v,\vec w(s))
$$
such that 
there are indices $i$ and $j$ with $b_{i1}(s)=1$, $b_{j1}(s)=-1$ and $b_{l1}=0$ if $l\ne i$ and $l\ne j$.
\end{Lemma}

\begin{pf} By Lemma \ref{Lemma5}, there exists a sequence of allowable column additions
$(A,\vec v,\vec w)\rightarrow (A(t_1),\vec v,\vec w(t_1))$
such that at most two entries of the first column of $B(t_1)$ are nonzero. Without loss of generality, we may assume that $b_{k1}=0$ if $k\ne 1$ or 2.

First assume that one of $b_{11}$ or $b_{21}$ is zero. We may suppose that $b_{21}=0$. Then since $\text{Det}(B)=\pm1$, we  have that $b_{11}=\pm1$.
As in Remark \ref{Remark2}, we can (if necessary) perform the permissible column addition $C_{21}$ a finite number of times so that the column
addition $C_{12}$ is permissible. We can then perform $C_{12}$ to get a matrix which satisfies the conclusions of the lemma in this case.

 Now assume that both  $b_{11}$ and $b_{21}$ are nonzero. 
 Since $b_{11}C_1+b_{21}C_2=e_1$, where  $e_1=(1,0,\ldots,0)^t$, it follows
 that $b_{11}$ and $b_{21}$ have opposite signs.
Recall that $\beta=\max\{|b_{11}|, |b_{21}|\}$.
 If $\beta=1$, then we have obtained the conclusions of the theorem. 
 
 Assume that $\beta>1$.  
We will show that we can construct a sequence of column additions in the first 3 columns which are *-allowable 
\begin{equation}\label{eq7}
(A,\vec v,\vec w)\rightarrow (A(1),\vec v,\vec w(1))\rightarrow \cdots
\rightarrow (A(s_1),\vec v,\vec w(s_1))
\end{equation}
 such that 
 $\beta(s_1)<\beta$.
 
 Once we have established the existence of the sequence (\ref{eq7}), we can apply Lemma \ref{Lemma5} to construct a sequence of allowable column additions 
 \begin{equation}\label{eq8}
(A(s_1),\vec v,\vec w(s_1))\rightarrow (A(s_1+1),\vec v,\vec w(s_1+1))\rightarrow \cdots
\rightarrow (A(s_2),\vec v,\vec w(s_2))
\end{equation}
 such that at most two of the entries in the first column of $B(s_2)$ are nonzero, and 
 $\beta(s_2)\le \beta(s_1)<\beta$.
 We  can thus alternate sequences (\ref{eq7}) and (\ref{eq8}) to eventually obtain the conclusions of the theorem.
 
 It remains to prove that we can construct a sequence (\ref{eq7}).

  Since $\text{Det}(B)=\pm 1$, and $\beta>1$,
 we must have that 
the maximum $\beta$ is obtained by  only one of $|b_{11}|$ and $|b_{21}|$. Without loss of generality, we may
 assume that 
 $$
 \beta=|b_{11}|> |b_{21}|.
 $$
 
  We now perform a finite sequence of *-allowable column additions  $C_{32}$, followed by a *-allowable column addition $C_{23}$ to obtain  a sequence of transformations of triples  
 $$
  (A,\vec v,\vec w)\rightarrow  \cdots
 \rightarrow (A(t_1),\vec v,\vec w(t_1))
 $$
 where the first column of $B(t_1)$ is
 $$
(b_{11}(t_1),b_{21}(t_2),\ldots,b_{n1}(t_1))^t=(b_{11},b_{21},-b_{21},0,\ldots,0)^t.
$$
with  $\beta (t_1) = |b_{11}(t_1)| = |b_{11}| = \beta $,
and  either $C_{13}$ or $C_{31}$ is allowable.  

If $C_{31}$ is allowable on $(A(t_1),\vec v,\vec w(t_1))$ we perform it to get 
$$
|b_{11}(t_1+1)| = |b_{11}(t_1)-b_{31}(t_1)| = |b_{11}+b_{21}|<\beta (t_1) = \beta 
$$
 and we stop.

If not, we have $C_{13}$ is allowable and after that, 
$b_{11}(t_1+1)= b_{11}$ and 
$$
b_{31}(t_1+1) = b_{31}(t_1)-b_{11}(t_1) = -b_{21}-b_{11}
$$
 have opposite signs. Further, $\beta (t_1+1) = \beta (t_1)$ and 
$w_3(t_1+1) = w_3(t_1)-w_1 \leq w_3-w_1$.  Now, $C_{32}$ or $C_{23}$ must be allowable.  

As before we perform a finite sequence of *-allowable column additions $C_{32}$, followed by a *-allowable column addition $C_{23}$ to obtain a sequence of transformations of triples  
 $$
  (A(t_1+1),\vec v,\vec w(t_1+1))\rightarrow  \cdots
 \rightarrow (A(t_2),\vec v,\vec w(t_2))
 $$
 where 
 $\beta(t_2) =\beta (t_1+1)= \beta$,
 $$
  \max\{\mid b_{21}(t_2)\mid,\mid b_{31}(t_2)\mid\}<\mid b_{11}(t_2)\mid=\beta(t_2)=\beta,
  $$
and one of 
 $C_{13}, C_{31}, C_{12}$ or  $C_{21}$ must now be allowable. 
  
  Performing an allowable $C_{31}$ or $C_{21}$ decreases $\beta$ and we stop. 
  If not, we perform $C_{13}$ or $C_{12}$ to get $\beta (t_2+1) = \beta (t_2)$ and none of the four 
  $C_{13}, C_{31}, C_{12}$ and $C_{21}$ are allowable.  
  
   Further, $w_2(t_2+1)$ or $ w_3(t_2+1)$ is reduced by  $w_1$, so that, 
$$
w_2(t_2+1)+w_3(t_2+1) = w_2(t_2)+w_3(t_2)-w_1 \leq w_2+w_3-2w_1.
$$
Now either $C_{32}$ or $C_{32}$ becomes allowable and we repeat this process.  
Since we can perform a $C_{13}$ or a $
C_{12}$ at most $[\frac{w_2+w_3}{w_1}]$ times, we must achieve a reduction in $\beta$ after a finite number of steps.

\end{pf}

\begin{Lemma}\label{Lemma10} Let $(A,\vec v,\vec w)$ be a triple such that $A=(C_1,\ldots, C_n)$ satisfies the relation
$$
C_k=C_1-e_1
$$
for some $k$, where $e_1=(1,0,\ldots,0)^t$.
 Let $A_{11}$ be the matrix obtained from $A$ be deleting the first row and column.
Then
$$
\text{Det}(A_{11})=\text{Det}(A)=\pm1.
$$ 
Let $\tilde v = (v_2,\ldots, v_n)^t$ and  $\tilde w=(\tilde w_2,\ldots,\tilde w_n)=A_{11}^{-1}\tilde v$. Then
$$
\tilde w_j=w_j\text{ for }j\ne k
$$
and 
$$
\tilde w_k=w_1+w_k.
$$
\end{Lemma}

\begin{pf} Set $\lambda=\text{Det}(A)=\pm1$. Subtracting the $k$-th column of $A$ from the first column, we see that $\text{Det}(A_{11})=\lambda$.
We thus have that $B=A^{-1}=\lambda\text{adj}(A)$, and $A_{11}^{-1}=\lambda\text{adj}(A_{11})$. Let
$$
A_{11}^{-1}=\lambda\text{adj}(A_{11})=\left(\begin{array}{llll}
x_{22}&x_{23}&\cdots&x_{2n}\\
x_{32}&x_{33}&\cdots&x_{3n}\\
&\vdots&&\\
x_{n2}&x_{n3}&\cdots&x_{nn}
\end{array}\right).
$$
Since $C_k=C_1-e_1$ and $\text{adj}(A)=\lambda A^{-1}$, the first column of $\text{adj}(A)$ is 
$$
(\lambda,0,\ldots,0,-\lambda,0,\ldots,0)^t,
$$
 where  $-\lambda$ occurs in the $k$-th row.

We will compute the entry $\lambda b_{ij}$ in the $i$-th row and $j$-th column of $\text{adj}(A)$. Let $A_{ji}$ be the matrix obtained from $A$ by deleting the $j$-th row and $i$-th column. 

First suppose that $i\ne 1$, $i\ne k$  and $j>1$. Subtracting the $k$-th column   of $A_{ji}$ (the $(k-1)$-st column of $A_{ji}$ if $i<k$)
from the first column, and expanding the determinant along the first column, we get that the $(i,j)$-th entry of $\text{adj}(A)$ is
$$
\begin{array}{ll}
(-1)^{i+j}\text{Det}(A_{ji})&=(-1)^{i+j}\text{Det}[(A_{ji})_{11}]\\
&=(-1)^{i+j}\text{Det}[(A_{11})_{j-1,i-1}]\\
&=\lambda x_{i,j}.
\end{array}
$$

For $j>1$, set
$$
t_j=\lambda(-1)^{1+j}\text{Det}(A_{j1}).
$$
We have that the $(1,j)$-th entry of $\text{adj}(A)$ is $\lambda t_j$.
We expand
$$
(-1)^{1+j}\text{Det}(A_{j1})=(-1)^{1+j}
\text{Det}\left(\begin{array}{llll}
a_{12}&\cdots&a_{11}-1&\cdots\\
&\vdots&&\\
a_{n2}&\cdots &a_{n1}&\cdots
\end{array}\right)
$$
$$
=(-1)^{1+j}
\text{Det}\left(\begin{array}{llll}
a_{12}&\cdots&a_{11}&\cdots\\
&\vdots&&\\
a_{n2}&\cdots &a_{n1}&\cdots
\end{array}\right)
+(-1)^{1+j+1+k-2}\text{Det}[(A_{11})_{j-1,k-1}]
$$
$$
=(-1)^{1+j+k-2}  \text{Det}(A_{jk})+ 
(-1)^{j+k}\text{Det}[(A_{11})_{j-1,k-1}]
$$
$$
=-(-1)^{j+k}\text{Det}(A_{jk})+\lambda x_{kj}.
$$
We see that, for $j>1$, the $(k,j)$-th entry  of $\text{adj}(A)$ is $\lambda (x_{k,j}-t_j)$.

In conclusion,
$$
B=A^{-1}=\lambda\text{adj}(A)=\left(\begin{array}{llll}
1&t_2&\cdots& t_n\\
0&x_{22}&\cdots&x_{2n}\\
&\vdots&&\\
-1&x_{k2}-t_2&\cdots&x_{kn}-t_n\\
&\vdots&&\\
0&x_{n2}&\cdots&x_{nn}
\end{array}\right).
$$
Now we see that
$$
\begin{array}{ll}
\tilde w_i&=(x_{i2},\ldots,x_{in})(v_2,\ldots,v_n)^t\\
&=(0,x_{i2},\ldots,x_{in})(v_1,\ldots,v_n)^t\\
&=w_i
\end{array}
$$
if $i\ne k$, and
$$
\begin{array}{ll}
\tilde w_k&=\sum_{j=2}^nx_{kj}v_j\\
&=[-v_1+\sum_{j=2}^n(x_{kj}-t_j)v_j]+[v_1+\sum_{j=2}^nt_jv_j]\\
&=w_k+w_1
\end{array}
$$
 \end{pf}

Now we prove Theorem \ref{Theorem1}. 

A quadruple $(A,\vec v,\vec w,k)$ is a triple $(A,\vec v,\vec w)$ and a number $k$ with $1<k\le n$ such that
if $A=(C_1,\ldots, C_n)$, then $C_k=C_1-e_1$.  To a quadruple $(A,\vec v,\vec w,k)$ we associate an $(n-1)$-dimensional triple $(\tilde A=A_{11},\tilde v,\tilde w)$ with the notation
of Lemma \ref{Lemma10}. 

A permissible transformation for the quadruple $(A,\vec v,\vec w,k)$ is a series of permissible column additions and row interchanges which transform the triple $(A,\vec v,\vec w)$
to a triple $(A(1),\vec v,\vec w(1))$ and a number $j$ with $1<j\le n$, such that $(A(1),\vec v,\vec w(1),j)$ is a quadruple.

By Lemma \ref{Lemma6}, and since $B=A^{-1}$, there exists a sequence of permissible column additions,  possibly followed by some row interchanges $T_{ij}$,
$(A,\vec v,\vec w)\rightarrow (A(1),\vec v,\vec w(1))$,
and a number $k(1)$ such that $(A(1),\vec v,\vec w(1),k(1))$ is a quadruple. Without loss of generality, we may assume that there exists a number $k$ such 
that $(A,\vec v,\vec w,k)$ is a quadruple.

If $n=3$, then after expanding the determinant of $\tilde A$, we see that after possibly performing the row interchange $T_{23}$,
there is a sequence of permissible row subtractions $R_{ji}$ which transform $\tilde A$ into the identity matrix.  
If $n>3$, we assume by induction that there exists 
 a sequence of permissible column additions and row interchanges 
\begin{equation}\label{eq11}
(\tilde A,\tilde v,\tilde w)\rightarrow (\tilde A(1),\tilde v, \tilde w(1)) \rightarrow \cdots \rightarrow (\tilde A(s), \tilde v,\tilde w (s))
\end{equation}
followed by a sequence of permissible row subtractions 
\begin{equation}\label{eq12}
(\tilde A(s),\tilde v, \tilde w(s))\rightarrow (\tilde A(s+1),\tilde v(s+1), \tilde w(s)) \rightarrow \cdots \rightarrow (\tilde A(l), \tilde v(l), \tilde w(s))
\end{equation}
such that $\tilde A(l)$ is the $(n-1)\times (n-1)$ identity matrix.

We will first construct a sequence of permissible transformations of quadruples 
\begin{equation}\label{eq13}
(A,\vec v,\vec w,k)\rightarrow (A(1),\vec v, \vec w(1), k(1)) \rightarrow \cdots \rightarrow (A(s), \vec v,\vec w (s), k(s))
\end{equation}
such that for $1\le t\le s$, we have that $A(t)_{11} = \tilde A(t)$ and $\tilde w(t)=A(t)_{11}^{-1}(v_2,\ldots,v_n)^t$.

Suppose that we have constructed  (\ref{eq13}) out to $(A(t),\vec v, \vec w(t), k(t))$, and $t<s$. We will construct 
$(A(t+1),\vec v, \vec w(t+1), k(t+1))$. 

First suppose that $\tilde A(t+1)$ is obtained from $\tilde A(t)$ by interchanging the $i$-th and $j$-th row. 
Let the triple $(A(t+1),\vec v(t+1),\vec w(t+1))$ be obtained from the triple   $(A(t),\vec v,\vec w(t))$ by performing the row interchange $T_{ij}$.
Then the row interchange $T_{ij}$
determines a permissible transformation of $(A(t),\vec v,\vec w(t),k(t))$ to $(A(t+1),\vec v,\vec w(t+1),k(t))$,
such that $A(t+1)_{11} = \tilde A(t+1)$.

Suppose that $\tilde A(t+1)$ is obtained from
$\tilde A(t)$ by adding the $j$-th column of $\tilde A(t)$ to the $i$-th column. We necessarily have that $\tilde w_j(t)>\tilde w_i(t)$. Set $k=k(t)$.

If $i\ne k$ and $j\ne k$ then we have (by Lemma \ref{Lemma10}) that $\vec w_j(t)>\vec w_i(t)$, and thus the column addition $C_{ij}$
determines a permissible transformation of $(A(t),\vec v,\vec w(t),k(t))$ to $(A(t+1),\vec v,\vec w(t+1),k(t))$,
such that $A(t+1)_{11} = \tilde A(t+1)$.

Suppose that $i=k$.
Then $\tilde w_j(t)>\tilde w_k(t)$. Since $\tilde w_k(t)= w_1(t)+w_{k}(t)$ and $\tilde w_j(t)= w_j(t)$
(by Lemma \ref{Lemma10}), we can construct a permissible transformation of
quadruples $(A(t),\vec v,\vec w(t),k(t))\rightarrow (A(t+1),\vec v,\vec w(t+1),k(t))$ by first performing the permissible column addition $C_{kj}$ followed
by the permissible column addition $C_{1j}$.
We have that $A(t+1)_{11} = \tilde A(t+1)$.

Suppose that $j=k$.  Then $\tilde w_k(t)>\tilde w_i(t)$. 

If $w_1(t)>w_i(t)$, then 
we define a permissible transformation of
quadruples 
$$
(A(t),\vec v,\vec w(t),k(t))\rightarrow (A(t+1),\vec v,\vec w(t+1),k(t))
$$
 by  performing the permissible column addition $C_{i1}$. 

Suppose that $w_1(t)<w_i(t)$. 
If ($\overline A,\vec v,\overline w)$ is the triple obtained from $(A(t),\vec v,\vec w(t))$ by $C_{1i}$, then we have that the $i$-th coefficient of $\overline w$ is $\overline w_i=w_i(t)-w_1(t)$. Since $\tilde w_k(t)>\tilde w_i(t)$, we must have that $w_1(t)+w_k(t)>w_i(t)$, which implies that 
$\overline w_k(t)>\overline w_i(t)$. 
Thus we can construct a permissible transformation of
quadruples 
$$
(A(t),\vec v,\vec w(t),k(t))\rightarrow (A(t+1),\vec v,\vec w(t+1),k(t+1)=i)
$$
 by first performing the permissible column addition $C_{1i}$ followed
by the permissible column addition $C_{ik}$.
We have that
$A(t+1)_{11} = \tilde A(t+1)$.

We can thus inductively construct the sequence (\ref{eq13}). 
Let $k=k(s)$.
Since $C_k(s)=C_1(s)-e_1$, where $C_k(s)$, $C_1(s)$ are the $k$-th and first columns of $A(s)$, the sequence of permissible row subtractions of (\ref{eq12}) gives a sequence of permissible row subtractions
$(A(s),\vec v,\vec w(s))\rightarrow (A(l),\vec v(l),\vec w(s))$ such that 
$A(l)$ is a matrix where  $A(l)_{11}$ is an identity matrix, $a_{1k}(l)=a_{11}(tl)-1$, $a_{k1}(l)=1$, and $a_{i1}(l)=0$  if $i\ne 1$ and $i\ne k$.

Now we perform the successive permissible row subtractions on $A(l)$ of subtracting $a_{11}(l)-1$ times the $k$-th row from the first row, subtracting
$a_{1i}(l)$ times the $i$-th row from the first row for all $i\ne k$, and finally subtracting the first row from the $k$-th row, to transform $A(l)$
into the identity matrix. This completes the proof of Theorem \ref{Theorem1}.

We say that a column subtraction is permissible on$(A,\vec v,\vec w)$ if it leaves the entries of $A$ nonnegative. If we subtract the $i$-th column from the
$j$-th column, then $\vec v$ is unchanged, but the coefficient $w_i$ of $\vec w$ is added to $w_j$. So as a corollary to Theorem \ref{Theorem1}, or
by a simple modification of the proof of Theorem \ref{Theorem1}, we obtain:

\begin{Theorem}\label{Theorem20}
Suppose that $(A,\vec v,\vec w)$ is a triple. Then there exists a sequence of permissible column additions and interchanges, followed by a sequence of
permissible column subtractions that transforms $A$ to the identity matrix.
\end{Theorem}

\section{Local  factorization of birational extensions}

Suppose that $R$ is a regular local ring with quotient field $K$, and $\nu$ is a valuation  of $K$, with valuation ring $V$, such that 
$V$ dominates $R$ ($R\subset V$ and the maximal ideal of $V$ intersects $R$ in its maximal ideal). 
A monoidal transform of $R$ along $\nu$ is a regular local ring $R(1)$  such that $R(1)=R[\frac{P}{f}]_m$, where $P$ is a regular prime of $R$, $f\in P$ is such that
$\nu(f)=\min\{\nu(g)\mid g\in P\}$, and $m=\{g\in R[\frac{P}{f}]\mid\nu(g)>0\}$. 
We have that $V$ dominates $R(1)$
and $R(1)$ dominates $R$.

\begin{Theorem}\label{Theorem15} Suppose that $k$ is  a field, $k[x_1,\ldots,x_n]$, $k[y_1,\ldots,y_n]$ are polynomial rings and 
there exists a matrix $(a_{ij})$ of nonnegative integers satisfying
$$
x_i=\prod_{j=1}^ny_j^{a_{ij}}
$$
for $1\le i\le n$ with $\text{Det}(a_{ij})=\pm 1$. Let $R=k[x_1,\ldots,x_n]_{(x_1,\ldots,x_n)}$ and $S=k[y_1,\ldots,y_n]_{(y_1,\ldots,y_n)}$.
Suppose that $\nu$ is a rank 1 valuation  of $k(y_1,\ldots,y_n)$ with valuation ring $V$ which dominates $S$,  such that $\nu(y_1),\ldots,\nu(y_n)$ are rationally
independent. Then there exists a commutative diagram 
\begin{equation}\label{eq20}
\begin{array}{lllll}
&&T&&\\
&\nearrow&&\nwarrow&\\
R&&\rightarrow&&S
\end{array}
\end{equation}
such that $T$ is a regular local ring dominated by $V$, and the northeast and northwest arrows are products of monoidal transforms along $\nu$.
\end{Theorem}

\begin{pf} Let $A=(a_{ij})$, $\vec v =(\nu(x_1),\ldots,\nu(x_n))^t$ and $\vec w=(\nu(y_1),\ldots,\nu(y_n))$.
By Theorem \ref{Theorem1}, there exists a sequence of permissible column additions and row interchanges
$$
(A,\vec v,\vec w)\rightarrow (A(1),\vec v, \vec w(1)) \rightarrow \cdots \rightarrow (A(s), \vec v,\vec w (s))
$$
followed by a sequence of permissible row subtractions 
$$
(A(s),\vec v, \vec w(s))\rightarrow (A(s+1),\vec v(s+1), \vec w(s)) \rightarrow \cdots \rightarrow (A(t), \vec v(t), \vec w(t))
$$
such that $A(t)$ is the $n\times n$ identity matrix.

We will construct a diagram (\ref{eq20}), in which the northwest arrow is a product of monodial transforms along $\nu$, 
\begin{equation}\label{eq21}
S\rightarrow S(1)\rightarrow \cdots \rightarrow S(s)=T
\end{equation}
and the northeast arrow is a product of monoidal transforms along $\nu$ 
\begin{equation}\label{eq22}
R\rightarrow R(1)\rightarrow \cdots\rightarrow R(t-s)=T.
\end{equation}

We inductively construct (\ref{eq21}), with a system of regular parameters $(y_1(l),\ldots,y_n(l))$ in each $S(l)$,  
so that $x_i=\prod_j y_j(l)^{a_{ij}(l)}$ for $1\le i\le n$,  $\vec v(l) =(\nu(x_1),\ldots,\nu(x_n))^t=\vec v$ and
$\vec w(l)=(\nu(y_1(l)),\ldots,\nu(y_n(l)))^t$ for $1\le l\le s$.

Suppose that $A(l+1)$ is obtained from $A(l)$ by the row interchange $T_{ij}$.
We define $S(l+1)$ to be $S(l)$, and we interchange the regular parameters $x_i$ and $x_j$ of $R$.

Suppose that $A(l+1)$ is obtained from $A(l)$ by the permissible column addition $C_{ij}$.
We define $S(l+1)$ to be  the local ring of the blow up of
the prime ideal $(y_i(l),y_j(l))$ which is dominated by $V$.   Since $\nu(y_j(l))>\nu(y_i(l))$, we have that 
$$
S(l+1)=S(l)[\frac{y_j(l)}{y_i(l)}]_{(y_1(l+1),\ldots,y_n(l+1))}
$$
where
$$
y_k(l+1)=\left\{\begin{array}{ll}
y_k(l)&\text{ if }k\ne j\\
\frac{y_j(l)}{y_i(l)}&\text{ if }k=j
\end{array}\right.
$$
are regular parameters  in $S(l+1)$. 

We now inductively construct (\ref{eq22}), with a system of regular parameters $(x_1(l),\ldots,x_n(l))$ in each $R(l)$,  
so that $x_i(l)=\prod_j y_j(s)^{a_{ij}(l+s)}$ for $1\le i\le n$,   $\vec v(l+s) =(\nu(x_1(l)),\ldots,\nu(x_n(l)))^t$ and
$\vec w(l+s)=(\nu(y_1(s)),\ldots,\nu(y_n(s)))^t=\vec w(s)$ for $1\le l\le t-s$.

Suppose that $A(l+1+s)$ is obtained from $A(l+s)$ by the permissible row subtraction $R_{ji}$.
We define $R(l+1)$ to be  the local ring of the blow up of
the prime ideal $(x_i(l),x_j(l))$ which is dominated by $V$.   Since $x_{i}(l)$ divides $x_j(l)$ in $T$, we have that 
$$
R(l+1)=R(l)[\frac{x_j(l)}{(x_i(l)}]_{(x_1(l+1),\ldots,x_n(l+1))}
$$
where
$$
x_k(l+1)=\left\{\begin{array}{ll}
x_k(l)&\text{ if }k\ne j\\
\frac{x_j(l)}{x_i(l)}&\text{ if }k=j
\end{array}\right.
$$
are regular parameters  in $R(l+1)$. 

Since $A(t)=Id$, we have that $x_i(t-s)=y_i(s)$ for $1\le i\le n$ and thus $T$ satisfies the conclusions of the theorem.
\end{pf}

\begin{Remark} If $S\rightarrow S(1)$ is a monoidal transform of a regular local ring $S$, then $S$ is called an inverse monoidal transform of $S(1)$
(Chapter 6 of \cite{C2}). With the notation of Theorem \ref{Theorem15}, a permissible column subtraction of $A=(a_{ij})$ induces an inverse monoidal
transform  $R\rightarrow S(1)\rightarrow S$ of $S$. We can use Theorem \ref{Theorem20}, instead of Theorem \ref{Theorem1}, to prove Theorem  \ref{Theorem15}.
Theorem \ref{Theorem20} proves the equivalent statement that a diagram (\ref{eq20}) can be constructed, where the northwest arrow is factored by
a sequence of inverse monoidal transforms from $T$ to $S$.
\end{Remark}

\begin{Theorem}\label{Theorem24}
Suppose that $R\subset S$ are regular local rings, essentially of finite type over a field $k$ of characteristic zero, with a common quotient field $K$,
such that $S$ dominates $R$.
Let $V$ be a valuation ring of $K$ which dominates $S$. Then there exists a regular local ring $T$, with quotient field $K$, such that $T$ dominates $S$,
$V$ dominates $T$, and the inclusions $R\rightarrow T$ and $S\rightarrow T$ can be factored by sequences of monoidal transforms 
\begin{equation}\label{eq24}
\begin{array}{lllll}
&&V&&\\
&&\uparrow&&\\
&&T&&\\
&\nearrow&&\nwarrow&\\
R&&\rightarrow&&S.
\end{array}
\end{equation}
\end{Theorem}

\begin{pf} 
Let $r=\text{rank}(V)$. We can perform monoidal transforms on $R$ and $S$ so that the assumptions of Theorem 5.5 \cite{C2} hold.
 By Theorem 5.5 \cite{C2}, there exists a commutative diagram of regular local rings
$$
\begin{array}{lll}
R'&\rightarrow&S'\\
\uparrow&&\uparrow\\
R&\rightarrow&S
\end{array}
$$
such that $R'$, $S'$ have respective regular parameters $(z_1,\ldots,z_n)$, $(w_1,\ldots,w_n)$ satisfying the conclusions of Theorem 5.5 \cite{C2}.
In particular, there exists a matrix $a_{ij}$ such that $z_i=\prod_jw_j^{a_{ij}}$ for $1\le i\le n$, where $A=(a_{ij})$ has the block form
$$
A=\left(\begin{array}{lllllll}
G_1&&&&0&&\\
&Id&&&&&\\
&&G_2&&&&\\
&&&Id&&&\\
0&&&&\cdots&&\\
&&&&&Id&\\
&&&&&&G_r
\end{array}\right).
$$
Here,
$G_i=(g_{jk}(i))$ is an $s_i\times s_i$ matrix of determinant $\pm1$, so that we have
$$
\begin{array}{ll}
z_{t_1+\cdots+t_{i-1}+1}&=w_{t_1+\cdots+t_{i-1}+1}^{g_{11}(i)}\cdots w_{t_1+\cdots+t_{i-1}+s_i}^{g_{1s_i}(i)}\\
&\vdots\\
z_{t_1+\cdots+t_{i-1}+s_{i}}&=w_{t_1+\cdots+t_{i-1}+1}^{g_{s_i1}(i)}\cdots w_{t_1+\cdots+t_{i-1}+s_i}^{g_{s_is_i}(i)}\\
\end{array}
$$
for $1\le i\le r$. We further have that $\nu(z_{t_1+\cdots+t_{i-1}+1}),\cdots,\nu(z_{t_1+\cdots+t_{i-1}+s_i})$ are rationally
independent, and if $V_i=V\cap k(z_{t_1+\cdots+t_{i-1}+1},\cdots,z_{t_1+\cdots+t_{i-1}+s_i})$, then $V_i$ has
rank 1 (and rational rank $s_{i}$).

Let 
$$
\overline R_i=k[z_{t_1+\cdots+t_{i-1}+1},\ldots,z_{t_1+\cdots+t_{i-1}+s_i}]
_{(z_{t_1+\cdots+t_{i-1}+1},\ldots,z_{t_1+\cdots+t_{i-1}+s_i})},
$$
$$
\overline S_i=k[w_{t_1+\cdots+t_{i-1}+1},\ldots,w_{t_1+\cdots+t_{i-1}+s_i}]
_{(w_{t_1+\cdots+t_{i-1}+1},\ldots,w_{t_1+\cdots+t_{i-1}+s_i})}.
$$
By Theorem \ref{Theorem15}, there exists a regular local ring $\overline T_i$ which is dominated by $V_i$ and a commutative diagram
$$
\begin{array}{lllll}
&&V_i&&\\
&&\uparrow&&\\
&&\overline T_i&&\\
&\nearrow&&\nwarrow&\\
\overline R_i&&\rightarrow&&\overline S_i
\end{array}
$$
such that the northeast and northwest arrows are products of monoidal transforms.
By performing the corresponding sequences of monoidal transforms on $R'$ and $S'$ for $1\le i\le r$, we obtain the conclusions of the theorem.
\end{pf}

\vskip.5truein
\noindent
Department of Mathematics

\noindent University of Missouri

\noindent Columbia, MO  65211

\end{document}